\documentclass[aop,dvips]{arximspdf}
\usepackage{graphics}

\doi{10.1214/009117906000001132}
\volume{35}
\issue{5}
\pubyear{2007}
\firstpage{1931}
\lastpage{1949}

\makeatletter
\newtheorem{lemma}{Lemma}
\newtheorem{theorem}{Theorem}
\newproclaim{definition}{Definition}
\newproclaim{rem}{Remark}
\newproclaim{example}{Example}
\newtheorem{proposition}{Proposition}

\newproclaim{condition}[theorem]{Condition}
\makeatother

\begin{document}
\begin{frontmatter}

\title{On superconvergence of sums of free random~variables}
\runtitle{On superconvergence of sums of free random variables}

\begin{aug}
\author[A]{\fnms{Vladislav} \snm{Kargin}\corref{}\ead[label=e1]{kargin@cims.nyu.edu}}
\runauthor{V. Kargin}
\affiliation{Courant Institute of Mathematical Sciences}
\address[A]{Courant Institute of Mathematical Sciences\\
251 Mercer Street\\
New York, New York 10012\\
USA\\
\printead{e1}} 
\end{aug}

\received{\smonth{4} \syear{2006}}
\revised{\smonth{7} \syear{2006}}

%
\begin{abstract}
This paper derives sufficient conditions for superconvergence of sums of
bounded free random variables and provides an estimate for the rate of
superconvergence.
\end{abstract}

%

\end{frontmatter}

\section{Introduction}\label{s1}

Free probability theory is an interesting generalization of classical
probability theory to a noncommutative setting. It was introduced in the
mid-1980's by Voiculescu \cite{voiculescu83,voiculescu86,voiculescu87}
as a tool for studying type $\mathrm{II}_{1}$ von Neumann
algebras. In many respects, free probability theory parallels classical
probability theory. There exist analogues of the central limit theorem
\cite{voiculescu86}, the law of large numbers \cite{bercovicipata96} and
the classification of infinitely divisible and stable laws %
\cite{bercovicivoiculescu92,bercovicipata99}. On the other
hand, certain features of free and classical probability theories differ
strikingly. Let $S_{n}=n^{-1/2}\sum_{i=1}^{n}X_{i}$, where $X_{i}$ are
identically distributed and free random variables. Then the law of $S_{n}$
approaches the limit law in a completely different manner than in the
classical case. To illustrate this, suppose that the support of $X_{i}$
is $[-1,1]$. Take a positive number $\alpha<1$. Then, in the classical
case, the probability of $\{ \vert S_{n}\vert>\alpha n\} $ is
exponentially small, but not zero. In contrast, in the noncommutative
case, the probability becomes identically zero for all sufficiently large $n$.
This mode of convergence is called \textit{superconvergence}
in \cite{bercovicivoiculescu95}.

In this paper, we extend the superconvergence result to a more general
setting of nonidentically distributed variables and estimate the rate of
superconvergence quantitatively. It turns out, in particular, that the
support of $S_{n}$ can deviate from the supporting interval of the limiting
law by not more than $c/\sqrt{n}$ and we explicitly estimate the
constant $c$. An example shows that the rate $n^{-1/2}$ in this estimate cannot be
improved.

Related results have been obtained in the
random matrix literature. For example,
\cite{johnstone01} considers the distribution of the largest
eigenvalue of
an empirical covariance matrix for a sample of Gaussian vectors. This
problem can be seen as a problem concerning the edge of the spectrum of a
sum of $n$ random rank-one operators in the $N$-dimensional vector space. More
precisely, the question concerns sums of the form
$S_{n}=\sum_{i=1}^{n}x_{i}x_{i}^{\prime}$, where $x_{i}$ is a random $N$-vector
with entries distributed according to the Gaussian law with the
normalized variance $1/N$. Then $S_{n}$ is a matrix-valued random variable
with the Wishart distribution.

Johnstone is interested in the asymptotic behavior of the distribution of
the largest eigenvalue of $S_{n}$. The asymptotics are derived under the
assumptions that both $n$ and $N$ approach $\infty$ and that
$\lim n/N=\gamma>0$, $\gamma\neq\infty$. Johnstone found that the largest
eigenvalue has variance of the order $n^{-2/3}$ and that after an
appropriate normalization, the distribution of the largest eigenvalue
approaches the Tracy--Widom law. This law has a right-tail asymptotically
equivalent to $\exp[ -( 2/3) s^{3/2}] $ and, in
particular, is unbounded from above. Johnstone's results have
generalized the
original breakthrough results in \cite{tracywidom96}
(see also \cite{tracywidom00}) for self-adjoint random matrices without
covariance structure. In \cite{soshnikov99} and \cite{soshnikov02},
it is
shown that the results regarding the asymptotic distribution of the largest
eigenvalue remain valid even if the matrix entries are not necessarily
Gaussian.

An earlier contribution, \cite{baisilverstein98}, also considered
empirical covariance matrices of large random vectors that are
not necessarily Gaussian and studied their largest eigenvalues. Again,
both $n $ and $N$ approach infinity and $\lim n/N=\gamma>0$,
$\gamma\neq \infty$. In contrast to
Johnstone, Bai and Silverstein were interested in the
behavior of the largest eigenvalue along a sequence of increasing random
covariance matrices. Suppose that the support of the limiting eigenvalue
distribution is contained in the interior of a closed interval, $I$.
Bai and
Silverstein showed that the probability that the largest eigenvalue lies
outside $I$ is zero for all sufficiently large $n$.

These results are not directly comparable with ours for several reasons.
First, in our case, the edge of the spectrum is not random in the classical
sense and so it does not make sense to talk about its variance. Second,
informally speaking, we are looking at the limit situation when
$N=\infty$, $n\rightarrow\infty$. Because of this, we use much easier techniques
than all of the aforementioned papers,
as we do not need to handle the interaction of the
randomness and the passage to the asymptotic limit. Despite these
differences, comparison of our results with the results of the random matrix
literature is stimulating. In particular, superconvergence in free
probability theory can be thought as an analogue of the Bai--Silverstein
result.

The rest of the paper is organized as follows. Section \ref{s2} provides the
necessary background about free probability theory and describes the main
result, Section~\ref{s3} recalls some results that we will need in the proof and
Section~\ref{s4} is devoted to the proof of the main result.

\section{Main theorem}\label{s2}
\begin{definition}
\emph{A noncommutative probability space} is a pair $( \mathcal{A},E) $,
where $\mathcal{A}$ is a unital $C^{\ast}$-algebra of bounded
linear operators acting on a complex separable Hilbert space and $E$ is a
linear functional from $\mathcal{A}$ to complex numbers. The operators
belonging to the algebra $\mathcal{A}$ are called \emph{noncommutative random
variables} or
simply \emph{random variables} and the functional $E$ is called the
\emph{expectation}.
\end{definition}

An algebra of bounded linear operators is a unital $C^{\ast}$-algebra
if it
contains the identity operator $I$ and if it is closed with respect to
the $\ast$-operation, that is, if \mbox{$A\in\mathcal{A}$},
then $A^{\ast}\in \mathcal{A}$,
where $A^{\ast}$ is the adjoint of operator $A$. The algebra is also
assumed to be closed with respect to convergence in the operator norm.
The definition of noncommutative random variables can be generalized to
include unbounded linear operators affiliated with algebra
$\mathcal{A}$;
for details, see \cite{bercovicivoiculescu93,maassen92}. In this
paper, we restrict our attention to bounded random variables.

The linear functional $E$ is assumed to satisfy the following
properties (in
addition to linearity): (i) $E(I)=1$; (ii) $E(A^{\ast})=E(A)$;
(iii) $E(AA^{\ast})\geq0$;
(iv)~$E(AB)=E(BA)$;
(v) $E(AA^{\ast})=0$ implies $A=0$;
and (vi) if $A_{n}\rightarrow A$, then
$E( A_{n}) \rightarrow E( A) $.

For each self-adjoint operator $A$, the expectation induces a continuous
linear functional on the space of continuous functions,
$E_{A}\dvtx f\rightarrow Ef( A) $, and by the Riesz theorem,
we can write this functional as a
Stieltjes' integral of $f$ over a measure. We call this measure, $\mu
$, the
\emph{measure associated with operator }$A$ \emph{and expectation }$E$.
If $%
P( d\lambda) $ is the spectral resolution associated with
operator $A$, then $\mu( d\lambda) =E( P(d\lambda) ) $.
It is easy to check that $\mu$ is a probability measure
on $\mathbb{R}$. If $A$ is a bounded operator and $\Vert A\Vert\leq L$, then
the support of $\mu$ is contained in the circle $\vert\lambda
\vert\leq L$.

The most important concept in free probability theory is that of free
independence of noncommuting random variables. Let a set of r.v.'s
$A_{1},\ldots,A_{n}$ be given. With each of them, we can associate an
algebra $\mathcal{A}_{i}$, which is generated by $A_{i}$; that is, it is the closure
of all polynomials in variables $A_{i}$ and $A_{i}^{\ast}$. Let
$\overline{A}_{i}$ denote an arbitrary element of algebra $\mathcal{A}_{i}$.
\begin{definition}
The algebras $\mathcal{A}_{1},\ldots,\mathcal{A}_{n}$ (and variables
$A_{1},\ldots,A_{n}$ that generate them) are said to be \emph{freely independent}
or \emph{free} if the following condition holds:
\[
\varphi\bigl( \overline{A}_{i(1)},\ldots,\overline{A}_{i(m)}\bigr) =0,
\]
provided that $\varphi( \overline{A}_{i(s)}) =0$ and
$i(s+1)\neq i(s)$.
\end{definition}

In classical probability theory, one of the most important theorems is the
central limit theorem (CLT). It has an analogue in noncommutative
probability theory.
\begin{proposition}
Let r.v.\textup{'}s $X_{i}$, $i=1,2,\ldots,$ be self-adjoint and free. Assume that
$E(X_{i})=0$, $\Vert X_{i}\Vert\leq L$ and
$\lim_{n\rightarrow\infty}[ E( X_{1}^{2}) +\cdots+E( X_{n}^{2})]/n=a_{2}$.
Then measures associated with r.v.\textup{'}s $n^{-1/2}\sum_{i=1}^{n}X_{i}$
converge in distribution to an absolutely continuous measure with density
\[
\phi(x)=\frac{1}{2\pi a_{2}}\sqrt{4a_{2}-x^{2}}
\chi_{\lbrack-2\sqrt{a_{2}},2\sqrt{a_{2}}]}( x) .
\]
\end{proposition}

This result was proven in \cite{voiculescu83} and later generalized in
\cite{maassen92} to unbounded identically distributed variables that
have a
finite second moment. Other generalizations can be found in \cite{pata96}
and \cite{voiculescu98}.

In the classical case, the behavior of large deviations from the CLT is
described by the Cram\'er theorem, the Bernstein inequality and their
generalizations. It turns out that in the noncommutative case, the behavior
of large deviations is considerably different. The theorem below gives some
quantitative bounds on how the distribution of a sum of free random
variables differs from the limiting distribution.

Let $X_{n,i}$, $i=1,\ldots,k_{n}$, be a double-indexed array of bounded
self-adjoint random variables. The elements of each row,
$X_{n,1},\ldots,X_{n,k_{n}}$, are assumed to be free, but are not necessarily
identically distributed. Their associated probability measures are
denoted $\mu_{n,i}$, their Cauchy transforms are $G_{n,i}( z) $,
their $k$th moments are~$a_{n,i}^{( k) }$, and so on. We define
$S_{n}=X_{n,1}+\cdots+X_{n,k_{n}}$ 
and the probability measure $\mu_n$ as the spectral probability measure
of $S_n$. We are interested in the behavior of probability measure
$\mu_n$ as $n$ grows.

We will assume that the first moments of the random variables $X_{n,i}$ are
zero and that $\Vert X_{n,i}\Vert\leq L_{n,i}$. Let
$v_{n}=a_{n,1}^{( 2) }+\cdots+a_{n,k_{n}}^{( 2) }$,
$L_{n}=\max_{i}\{ L_{n,i}\} $ and
$T_{n}=L_{n,1}^{3}+\cdots+L_{n,k_{n}}^{3}$.
\begin{theorem}\label{main_theorem}
Suppose that $\lim\sup_{n\rightarrow\infty}T_{n}/v_{n}^{3/2}<2^{-12}$.
Then for all sufficiently large $n$, the support of $\mu_{n}$ belongs to
\[
I=\bigl( -2\sqrt{v_{n}}-cT_{n}/v_{n},2\sqrt{v_{n}}+cT_{n}/v_{n}\bigr) ,
\]
where $c>0$ is an absolute constant.
\end{theorem}
\begin{rem}\label{rem1}
$c=5$ will suffices although it is not the best possible.
\end{rem}
\begin{rem}\label{rem2}
Informally, the assumption that
$\lim\sup_{n\rightarrow\infty}T_{n}/v_{n}^{3/2}<2^{-12}$
means that there are no large outliers. An
example of when the assumption is violated is provided by random variables with
variance $a_{n,i}^{( 2) }=n^{-1}$ and $L_{n,i}=1$. Then $T_{n}=n$
and $v_{n}^{3/2}=1$, so that $T_{n}/v_{n}^{3/2}$ increases when $n$ grows.
\end{rem}
\begin{rem}\label{rem3}
The assumptions in Theorem \ref{main_theorem} are weaker than the
assumptions in Theorem 7 of \cite{bercovicivoiculescu95}. In particular,
Theorem \ref{main_theorem} allows us to draw conclusions about random variables
with nonuniformly bounded support. Consider, for example, random
variables $X_{k}$, $k=1,\ldots,n$, that are supported on intervals
$[-k^{1/3},k^{1/3}] $ and have variances of order $k^{2/3}$.
Then $T_{n}$ has the order
of $n^{2}$ and $v_{n}$ has the order of~$n^{5/3}$. Therefore,
$T_{n}/v_{n}^{3/2}$ has the order of $n^{-1/2}$ and Theorem \ref{main_theorem}
is applicable. It allows us to conclude that the support of
$S_{n}=X_{1}+\cdots+X_{n}$ is contained in the interval
$( -2\sqrt{v_{n}}-cn^{1/3}, 2\sqrt{v_{n}}+cn^{1/3})$.
\end{rem}
\begin{example}[(\textit{Identically distributed variables})]
A particular case of the above scheme involves
the normalized sums of identically
distributed, bounded, free r.v.'s
$S_{n}=( \xi_{1}+\cdots+\xi_{n}) /\sqrt{n}$.
If $\Vert\xi_{i}\Vert\leq L$, then
$\Vert \xi _{i}/\sqrt{n}\Vert\leq L_{n,i}=L_{n}=L/\sqrt{n}$. Therefore,
$T_{n}=L^{3}/\sqrt{n}$. If the second moment of $\xi_{i}$ is
$\sigma^{2}$, then the second moment of the sum $S_{n}$ is $v_{n}=\sigma^{2}$. Applying
the theorem, we obtain the result that starting with certain $n$, the support
of the distribution of $S_{n}$ belongs to
$( -2\sigma-c(L^{3}/\sigma^{2}) n^{-1/2},2\sigma+c( L^{3}/\sigma^{2})n^{-1/2}) $.
\end{example}
\begin{example}[(\textit{Free Poisson})]
Let the $n$th row of our scheme have $k_{n}=n$ identically distributed
random variables $X_{n,i}$ with the Bernoulli distribution that places
probability $p_{n,i}$ on $1$ and $q_{n,i}=1-p_{n,i}$ on $0$. (It is
easy to
normalize this distribution to have the zero mean by subtracting $p_{n,i}$.)
Suppose that $\max_{i}p_{n,i}\rightarrow0$ as $n\rightarrow\infty$ and
that
\[
\sum_{i=1}^{n}p_{n,i}\rightarrow\lambda>0
\]
as $n\rightarrow\infty$. Then $L_{n,i}\sim1$ and $a_{n,i}^{2}=p_{n,i}
( 1-p_{n,i}) $ so that $T_{n}\sim n$ and $v_{n}\rightarrow\lambda$
as $n\rightarrow\infty$. Therefore, Theorem \ref{main_theorem}
does not apply. An easy calculation for the case $p_{n,i}=\lambda/n$ shows
that superconvergence still holds. This example shows that the
conditions of
the theorem are not necessary for superconvergence to hold.
\end{example}
\begin{example}[(\textit{Identically distributed binomial variables})]
Let $X_{i}$ be identically distributed with a distribution that
attributes positive weights $p$ and $q$ to $-\sqrt{q/p}$ and $\sqrt{p/q}$,
respectively. Then $EX_{i}=0$ and $EX_{i}^{2}=1$. It is not difficult to
show that the support of $S_{n}=n^{-1/2}\sum_{i=1}^{n}X_{i}$ is the interval
$I=[ x_{1},x_{2}] $, where
\[
x_{1,2}=\pm2\sqrt{1-\frac{1}{n}}+\frac{q-p}{\sqrt{pq}}\frac{1}{\sqrt{n}}.
\]
This example shows that the rate of $n^{-1/2}$ in Theorem \ref{main_theorem}
cannot be improved without further restrictions. Note, also, that for
$p>q$, $L_{n}$ is $\sqrt{p/q}$ and therefore the coefficient
preceding $n^{-1/2}$ is
of order $L_{n}$. In the general bound, the coefficient
is~$L_{n}^{3}/\sigma^{2}$. It is not clear whether it is possible to replace the
coefficient in
the general bound by a term of order $L_{n}$.
\end{example}

\section{Preliminary results}\label{s3}
\begin{definition}
The function
%
\begin{equation}\label{G_integral}
G( z) =\int_{-\infty}^{\infty}\frac{\mu( dt) }{z-t}
\end{equation}
is called the \emph{Cauchy transform} of the probability measure $\mu(dt) $.
\end{definition}

The Cauchy transform encodes a wealth of information about the underlying
probability measure. For our purposes, we need only some of them. First,
the following inversion formula holds.
\begin{proposition}[(The Stieltjes--Perron inversion formula)]
For any interval $[ a,b]$,
\[
\mu[ a,b] =-\lim_{\varepsilon\downarrow0}\frac{1}{\pi}
\int_{a}^{b}\operatorname{Im}G( x+i\varepsilon) \,dx,
\]
provided that $\mu( a) =\mu( b) =0$.
\end{proposition}

A proof can be found in \cite{akhiezer65}, pages 124--125.

We will call a function \emph{holomorphic} at a point $z$ if it can be
represented by a convergent power series in a sufficiently small disc with
center $z$. We call the function \emph{holomorphic in an open
domain} $D$ if it is holomorphic at every point of the domain. Here, $D$ may
include $\{ \infty\} $, in which case it is a part of the
extended complex plane $\mathbb{C}\cup\{ \infty\} $ with the
topology induced by the stereographic projection of the Riemann sphere
onto the extended complex plane.

The integral representation (\ref{G_integral}) shows that
the Cauchy transform
of every probability measure, $G( z) $, is a holomorphic function
in $C^{+}=\{ z\in\mathbb{C}|\operatorname{Im}z>0\} $ and
$C^{-}=\{z\in\mathbb{C}|\operatorname{Im}z<0\} $. If, in addition, the measure is
assumed to be supported on interval $[ -L,L] $, then the Cauchy
transform is holomorphic in the area $\Omega\dvtx\vert z\vert>L$
where it can be represented by a convergent power series in $z^{-1}$,
%
\begin{equation}\label{G_power_series}
G(z)=\frac{1}{z}+\frac{m_{1}}{z^{2}}+\frac{m_{2}}{z^{3}}+\cdots.
\end{equation}
Here, $m_{k}$ denote the moments of the measure $\mu$:
\[
m_{k}=\int_{-\infty}^{\infty}t^{k}\mu( dt) .
\]
In particular, $G( z) $ is holomorphic at
$\{ \infty\} $. We call series (\ref{G_power_series}) the $G$-series.

In the other direction, we have the following result.
\begin{lemma}\label{G_boundedness_lemma}
Suppose that:
\begin{longlist}[(2)]
\item[(1)] $G( z) $ is the Cauchy transform of a compactly supported
probability distribution, $\mu$, and
\item[(2)] $G( z) $ is holomorphic at every $z\in\mathbb{R}$,
$\vert z\vert>L$.
\end{longlist}
Then the support of $\mu$ lies entirely in the interval $[ -L,L] $.
\end{lemma}
\begin{pf}
From assumption (1), we infer that in some neighborhood of
infinity, $G( z) $ can be represented by the convergent power
series (\ref{G_power_series}) and that $G(z)$ is also holomorphic everywhere
in $C^{+}$ and $C^{-}$. Therefore, using assumption (2), we can conclude
that $G( z) $ is holomorphic everywhere in the area
$\Omega=\{ z|\vert z\vert>L\} $ including the point at infinity.

Let us detail the proof of this statement. Define%
\[
a=\inf\{ l\geq0|G( z) \mbox{ is holomorphic on }\vert z\vert>l\}
\]
and suppose, by seeking a contradiction, that $a>L$. Let $\varepsilon$ be
such that \mbox{$a-\varepsilon>L$}. Consider the area
$\Omega_{\varepsilon}=\{ z|a-\varepsilon<\vert z\vert <a+\varepsilon\} $.
Since $G(z)$ is holomorphic everywhere in $C^{+}$
and $C^{-}$, it is holomorphic in $\Omega_{\varepsilon}\setminus \mathbb{R}$.
In addition, by assumption (2), $G( z) $ is holomorphic at each
point of $\Omega_{\varepsilon}\cap\mathbb{R}$. Therefore, it is
holomorphic everywhere in $\Omega_{\varepsilon}$ and thus it is
holomorphic everywhere in
$\Omega_{\varepsilon}\cup\{ z|\vert z\vert>a+\varepsilon/2\}
=\{ z|\vert z\vert>a-\varepsilon\} $. This
contradicts the definition of $a$. Therefore, $a\leq L$ and
$G( z) $ is holomorphic everywhere in the area $\Omega= \{ z|\vert z\vert>L\} $,
including the point at infinity.

It follows that the power series (\ref{G_power_series}) converges everywhere
in the area $\Omega=\{ z|\vert z\vert
>L\} $. Since this power series has real coefficients, we can conclude
that $G(z)$ is real for $z\in R$, $\vert z\vert>L$. Also,
since $G(z)$ is holomorphic, and therefore continuous,
in $\vert z \vert >L, $ we can conclude that
$\lim_{\varepsilon\downarrow0}\operatorname{Im}G(z+i\varepsilon)=0$.
Then the Stieltjes inversion formula implies
$\mu( \lbrack a,b]) =0$ for each $[ a,b] \subset\{ x\in\mathbb{R}|\vert
x \vert>L\} $ provided that $\mu( a) =0$ and $\mu(b)=0$. It
remains to prove that this implies $\mu\{ \vert x\vert>L\} =0$.

For this purpose, note that the set of points $x\in\mathbb{R}$ for
which $\mu(x)>0$ is at most countable. Indeed, let $S$ be the set of all
$x$ for
which $\mu(x)>0$. We can divide this set into a countable collection of
disjoint subsets $S_{k}$, where $k$ are all positive integers and
$S_{k}=\{ x|k^{-1}\geq\mu\{ x\} >( k+1)^{-1}\} $.
Clearly, every $S_{k}$ is either empty or a finite set.
Otherwise, we could take an infinite countable sequence of
\mbox{$x_{i,k}\in S_{k}$} and would get (by countable additivity and monotonicity
of~$\mu$) that $\mu( S_{k}) \geq\sum_{i}\mu( x_{i,k})=+\infty$.
By monotonicity of $\mu$, we would further get
$\mu(\mathbb{R}) =+\infty$, which would contradict the assumption
that $\mu$ is a probability measure. Therefore, $S$ is a countable union of
finite sets $S_{k}$ and hence countable.

From the countability of $S$, we conclude that the set of points $x$ for
which $\mu( x) =0$ (i.e., $S^{c}$) is dense in the set $\vert
x\vert>L$. Indeed, take an arbitrary nonempty interval
$(\alpha,\beta) $. Then $( \alpha,\beta) \cap S^{c}\neq \varnothing$
since, otherwise, $( \alpha,\beta) \subset S$ and
therefore~$S$ would be uncountable. Hence, $S^{c}$ is dense. Using the denseness
of $S^{c}$, we can cover the set $\{ \vert x\vert>L\}$
by a countable union of disjoint intervals $[a,b]$, where
$\mu(a) =0$ and $\mu(b)=0$. For each of these intervals,
$\mu(\lbrack a,b]) =0$ and therefore countable additivity implies that
$\mu( \{ \vert x\vert>L\} ) =0$. Consequently,
$\mu$ is supported on a set that lies entirely in $[ -L,L] $.
\end{pf}
\begin{definition}
The inverse of the $G$-series (\ref{G_power_series})
(in the sense of formal series)
always exists and is called the \emph{K-series},
\[
G( K(z)) =K(G(z))=z.
\]
\end{definition}

In case of a bounded self-adjoint random variable $A$, the $G$-series are
convergent for $\vert z\vert\geq\Vert A\Vert$ and
the limit coincides with the Cauchy transform $G( z) $. As a
consequence, the $K$-series is convergent in a sufficiently small punctured
neighborhood of 0. We will call the limit $K( z) $. This function
has a pole of order $1$ at~$0$.

It is sometimes useful to know how functions $G( z) $ and
$K( z) $ behave under a rescaling of the random variable.
\begin{lemma}\label{dilation_properties}
\textup{(i)} $G_{\alpha A}(z)=\alpha^{-1}G_{A}(z/a)$ and
\textup{(ii)} $K_{\alpha A}(u)=\alpha K_{A}(\alpha u)$.
\end{lemma}

The claim of the lemma follows directly from definitions.

The importance of $K$-functions is that they allow us to compute the
distribution of the sum of free random variables.
\begin{proposition}[(Voiculescu's addition formula)]\label{additivity}
Suppose that self-adjoint r.v.\textup{'}s $A$ and $B$ are free. Let
$K_{A}$, $K_{B}$ and $K_{A+B}$ be the K-series for variables $A$, $B$
and $A+B$, respectively. Then
\[
K_{A+B}(u)=K_{A}(u)+K_{B}(u)-\frac{1}{u},
\]
where the equality holds in the sense of formal power series.
\end{proposition}

The proof can be found in \cite{voiculescu86}. Using this property, we can
compute the distribution of the sum of free r.v.'s as follows. Given two
r.v.'s, $A$ and $B$, compute their $G$-series. Invert them to obtain the $K$-series.
Use Proposition \ref{additivity} to compute $K_{A+B}$ and invert it to
obtain $G_{A+B}$. Use the Stieltjes inversion formula to compute the measure
corresponding to this $G$-series. This is the probability measure
corresponding to $A+B$.

\section[Proof of Theorem 1]{Proof of Theorem \protect\ref{main_theorem}}\label{s4}

The key ideas of the proof are as follows.
\begin{longlist}[(2)]
\item[(1)] We know that the Cauchy transform of the sum $S_{n}$ is the Cauchy
transform of a bounded r.v. (since, by assumption, each $X_{n,i}$ is bounded).
Consequently, the Cauchy transform of $S_{n}$ is holomorphic in a certain
circle around infinity (i.e., in the area $|z|>R$ for some $R>0$). We want
to estimate $R$ and apply Lemma \ref{G_boundedness_lemma} to conclude
that $S_{n}$ is supported on $[ -R,R] $.
\item[(2)] Since the $K$-function of $S_{n}$, call it $K_{n}(z)$, is the sum of
the $K$-functions of $X_{n,i}$ and the latter are functional inverses of the Cauchy
transforms of $X_{n,i}$, it is an exercise in complex analysis to prove that
the $K$-function of $S_{n}$ takes real values and is a one-to-one
function on
a sufficiently large real interval around zero. Therefore, it has a
differentiable functional inverse defined on a sufficiently large real
interval around infinity (i.e., on the set $I=( -\infty,-A]
\cup[ A,\infty) $ for some $A$ which we can explicitly estimate).
Moreover, with a little bit more effort, we can show that this inverse
function is well defined and holomorphic in an open complex
neighborhood of $I$. This shows that Lemma \ref{G_boundedness_lemma} is applicable, and the
estimate for $A$ provides the desired estimate for the support of $S_{n}$.
\end{longlist}

We will begin by finding the radius of convergence of the Taylor
series of
$K_{n}( z) $. First, we need to prove some preliminary facts about
Cauchy transforms of $X_{n,i}$.

Define $g_{n,i}( z) =G_{n,i}( z^{-1}) $. Since the
series $G_{n,i}( z) $ are convergent everywhere in
$\vert z\vert>L_{n,i}$, the Taylor series for $g_{n,i}(z) $
converges everywhere in $\vert z\vert<L_{n,i}^{-1}$.

Assume that $R_{n,i}$ and $m_{n,i}$ are such that:
\begin{enumerate}
\item$R_{n,i}\geq L_{n,i}$;
\item$\vert G_{n,i}( z) \vert\geq m_{n,i}>0$
everywhere on $\vert z\vert=R_{n,i}$;
\item$g_{n,i}( z) $ has only one zero in $\vert
z\vert<R_{n,i}^{-1}$.
\end{enumerate}

For example, we can take $R_{n,i}=2L_{n,i}$ and $m_{n,i}=(
4L_{n,i}) ^{-1}$. Indeed, for any $z$ with $| z| =r>L_{n,i}$,
we can estimate $G_{n,i}( z) $:
\begin{eqnarray*}
| G_{n,i}( z) | &\geq&\frac{1}{r}
- \biggl( \frac{a_{n,i}^{2}}{r^{3}}+\frac{| a_{n,i}^{3}| }{r^{4}}+\cdots\biggr)
\\
&\geq&\frac{1}{r}-\biggl( \frac{L_{n,i}^{2}}{r^{3}}
+ \frac{L_{n,i}^{3}}{r^{4}}+\cdots \biggr)
\\
&= &\frac{1}{r}-\frac{L_{n,i}^{2}}{r^{2}}\frac{1}{r-L_{n,i}}.
\end{eqnarray*}
In particular, taking $r=2L_{n,i}$, we obtain the estimate:
\[
| G_{n,i}( z) | \geq\frac{1}{4L_{n,i}},
\]
valid for every $i$ and everywhere on $| z| =2L_{n,i}$.

It remains to show that $g_{n,i}(z)$ has only one zero in $\vert
z\vert<( 2L_{n,i}) ^{-1}$. This is indeed so because
\[
g_{n,i}(z)= z\bigl( 1+a_{n,i}^{( 2) }z^{2}+a_{n,i}^{(3) }
z^{3}+ \cdots \bigr)
\]
and we can estimate
\[
\bigl\vert a_{n,i}^{( 2) }z^{2}+a_{n,i}^{( 3)}z^{3}+\cdots\bigr\vert
\leq L_{n,i}^{2}\biggl( \frac{1}{2L_{n,i}}\biggr)^{2}
+ L_{n,i}^{3}\biggl( \frac{1}{2L_{n,i}}\biggr) ^{3}+\cdots
= \frac{1}{2}.
\]

Therefore, Rouch\'e's theorem is applicable and $g_{n,i}$ has only one
zero in $\vert z\vert<( 2L_{n,i}) ^{-1}$.
\begin{definition}
Let $R_{n}=\max_{i}\{ R_{n,i}\} $, $m_{n}=\min_{i}\{m_{n,i}\} $
and $D_{n}=\break \sum_{i=1}^{k_{n}}R_{n,i}(m_{n,i})^{-2}$.
\end{definition}

We are now able to investigate the region of convergence for the series
$K_{n,i}( z) $. First, we prove a modification of Lagrange's
inversion formula.
\begin{lemma}\label{Lagranges_series}
Suppose $w=G(z)$ (where $G$ is not necessarily a
Cauchy transform) is holomorphic in a neighborhood of $z_{0}=\infty$ and
has the expansion
\[
G(z)=\frac{1}{z}+\frac{a_{1}}{z^{2}}+\cdots,
\]
converging for all sufficiently large $z$. Define $g(z)=G(1/z)$. Then the
inverse of $G( z) $ is well defined in a neighborhood of $0$ and
its Laurent series at $0$ is given by the formula
\[
z=G^{-1}( w) =\frac{1}{w}+a_{1}-\sum_{n=1}^{\infty}
\biggl[\frac{1}{2\pi in}
\oint_{\partial\gamma}\frac{dz}{z^{2}g(z)^{n}}\biggr] w^{n},
\]
where $\gamma$ is a sufficiently small disc around $0$.
\end{lemma}
\begin{pf}
Let $\gamma$ be a closed disc around $z=0$ in which $g(z)$
has only one zero. This disc exists because $g( z) $ is
holomorphic in a neighborhood of $0$ and has a nonzero derivative at $0$.
Let
\[
r_{w}=\tfrac{1}{2}\inf_{z\in\partial\gamma}\vert g( z)\vert.
\]
Then $r_{w}>0$, by our assumption on $\gamma$. We can apply Rouch\'e's theorem
and conclude that the equation $g( z) -w=0$ has only one solution
inside $\gamma$ if $\vert w\vert\leq r_{w}$. Let us consider a
$w$ such that $\vert w\vert\leq r_{w}$. Inside $\gamma$, the function
\[
\frac{g^{\prime}(z)}{z( g(z)-w) }
\]
has a pole at $z=1/G^{-1}(w)$ with the residue $G^{-1}(w)$ and a pole
at $%
z=0$ with the residue $-1/w$. Consequently, we can write:%
\[
G^{-1}( w) =\frac{1}{2\pi i}\oint_{\partial\gamma}
\frac{g^{\prime}(z)\, dz}{z( g(z)-w) }+\frac{1}{w}.
\]
The integral can be rewritten as follows:
\begin{eqnarray*}
\oint_{\partial\gamma}\frac{g^{\prime}(z)\,dz}{z( g(z)-w) }
&=&\oint_{\partial\gamma}\frac{g^{\prime}(z)}{zg(z)}
\frac{1}{1-{w}/{g( z) }}\,dz
\\
&=&\sum_{n=0}^{\infty}\oint_{\partial\gamma}
\frac{g^{\prime}(z)\,dz}{zg(z)^{n+1}}w^{n}.
\end{eqnarray*}
For $n=0$, we calculate
\[
\frac{1}{2\pi i}\oint_{\partial\gamma}\frac{g^{\prime}(z)\,dz}{zg(z)}
=a_{1}.
\]
Indeed, the only pole of the integrand is at $z=0$, of order two, and
the corresponding residue can be computed from the series expansion
for $g(z)$:
\begin{eqnarray*}
\mathrm{res}_{z=0}\frac{g^{\prime}(z)\,dz}{zg(z)}
&=& \frac{d}{dz}\frac{z^{2}( 1+2a_{1}z+\cdots) }
{z( z+a_{1}z^{2}+\cdots) }\bigg|_{z=0}
\\
&=& \frac{d}{dz}\frac{1+2a_{1}z+\cdots}{1+a_{1}z+\cdots}
\bigg|_{z=0}=a_{1}.
\end{eqnarray*}
For $n>0$, we integrate by parts:
\[
\frac{1}{2\pi i}\oint_{\partial\gamma}\frac{g^{\prime}(z)\,dz}
{zg(z)^{n+1}}
= -\frac{1}{2\pi i}\frac{1}{n}\oint_{\partial\gamma}
\frac{dz}{z^{2}g(z)^{n}}.
\]\upqed
\end{pf}
\begin{lemma}
The radius of convergence of the $K$-series for measure $\mu_{n}$ is at
least $m_{n}$.
\end{lemma}

The lemma essentially says that if r.v.'s $X_{n,1},\ldots,X_{n,k_{n}}$ are
all bounded by $L_{n}$, then the $K$-series for $\sum_{i}X_{n,i}$ converges
in the circle $\vert z\vert\leq1/( 4L_{n}) $.
\begin{pf}
Let us apply Lemma \ref{Lagranges_series} to
$G_{n,i}(z) $ with $\gamma$ having radius $( R_{n,i}) ^{-1}$. By
Lemma \ref{Lagranges_series}, the coefficients in the series for the inverse
of $G_{n,i}( z) $ are%
\[
b_{n,i}^{(k)}=\frac{1}{2\pi ik}\oint_{\partial\gamma}
\frac{dz}{z^{2}g_{n,i}(z)^{k}}
\]
and we can estimate them as
\[
\bigl\vert b_{n,i}^{( k) }\bigr\vert
\leq\frac{R_{n,i}}{k}( m_{n,i}) ^{-k}.
\]
This implies that the radius of convergence of the
$K$-series for measure $\mu_{n,i}$ is~$m_{n,i}$.
Consequently, the radius of convergence of the $K$-series
for measure $\mu_{n}$ is at least $m_{n}$.
\end{pf}

We can now investigate the behavior of $K_{n}( z) $ and its
derivative inside its circle of convergence.
\begin{lemma}\label{Lemma_K_estimate}
For every $z$ in $\vert z\vert<m_{n}$,
the following inequalities are valid:
%
\begin{eqnarray}
\biggl\vert K_{n}( z) -\frac{1}{z}-v_{n}z\biggr\vert
&\leq &D_{n}\vert z\vert^{2},
\label{K_estimate} \\
\biggl\vert K_{n}^{\prime}( z) +\frac{1}{z^{2}}-v_{n}\biggr\vert
&\leq&2D_{n}\vert z\vert.
\label{K_prime_estimate}
\end{eqnarray}
\end{lemma}

Note that $D_{n}$ is approximately $k_{n}L_{n}^{3}$, so the meaning of the
lemma is that the growth of $K_{n}-z^{-1}-v_{n}z$ around $z=0$ is
bounded by a constant that depends on the norm of the
variables $X_{n,1},\ldots,X_{n,k_{n}}$.
\begin{pf}
Consider the circle of radius $m_{n,i}/2$. We can estimate
$K_{n,i}$ inside this circle:
\begin{eqnarray*}
\biggl\vert K_{n,i}-\frac{1}{z}-a_{n,i}^{( 2) }z\biggr\vert
&\leq& \frac{R_{n,i}}{2}( m_{n,i}) ^{-2}\vert z\vert^{2}
+ \frac{R_{n,i}}{3}( m_{n,i}) ^{-3}\vert z\vert^{2}\frac{m_{n,i}}{2}
\\
&&{} +\frac{R_{n,i}}{4}( m_{n,i}) ^{-3}\vert z\vert^{2}
\frac{m_{n,i}^{2}}{2^{2}}+\cdots
\\
&=& R_{n,i}( m_{n,i}) ^{-2}\vert z\vert^{2}
\biggl(\frac{1}{2}+\frac{1}{3}\frac{1}{2}+\frac{1}{4}\frac{1}{2^{2}}
+\cdots\biggr)
\\
&\leq & R_{n,i}( m_{n,i}) ^{-2}\vert z\vert^{2}.
\end{eqnarray*}
Consequently, using Voiculescu's addition formula, we can estimate
%
\begin{equation}\label{K_estimate2}
\biggl\vert K_{n}( z) -\frac{1}{z}-v_{n}z\biggr\vert
\leq D_{n}\vert z\vert^{2}
\end{equation}
and a similar argument leads to the estimate
%
\begin{equation}\label{K_prime_estimate2}
\biggl\vert K_{n}^{\prime}( z) +\frac{1}{z^{2}}-v_{n} \biggr\vert
\leq2D_{n}\vert z\vert.
\end{equation}\upqed
\end{pf}
\begin{lemma}\label{K_zeros}
Suppose that \textup{(i)} $m_{n}>4/\sqrt{v_{n}}$ and
\textup{(ii)} $r_{n}\geq 4D_{n}/v_{n}^{2}$.
Then there are no zeros of $K_{n}^{\prime}(z) $
inside $\vert z\vert<1/\sqrt{v_{n}}-r_{n}$.
\end{lemma}

In other words, $K_{n}( z) $ has no critical points in a circle
which is sufficiently separated from $z=\pm1/\sqrt{v_{n}}$.
\begin{pf*}{Proof of Lemma \ref{K_zeros}}
If $r_{n}\geq v_{n}^{-1/2}$, then the set
$\vert z\vert<1/\sqrt{v_{n}}-r_{n}$ is empty and we are done. In the
following, we assume that $r_{n}<v_{n}^{-1/2}$. On
$\vert z \vert =v_{n}^{-1/2}-r_{n}$, we have \mbox{$\vert z\vert^{-2}>v_{n}$}.
Also, $\vert z-v_{n}^{-1/2}\vert\vert z+v_{n}^{-1/2}\vert
>r_{n}v_{n}^{-1/2}$. This is easy to see by considering the two cases
$\operatorname{Re}z\geq0$ and $\operatorname{Re}z\leq0$. In the first case,
$\vert z-v_{n}^{-1/2}\vert\geq r_{n}$ and
$\vert z+v_{n}^{-1/2}\vert>v_{n}^{-1/2}$. In the second case,
$\vert z-v_{n}^{-1/2}\vert>v_{n}^{-1/2}$ and
\mbox{$\vert z+v_{n}^{-1/2}\vert\geq r_{n}$}. Hence, in both cases,
the product $\vert z-v_{n}^{-1/2}\vert\vert z+v_{n}^{-1/2}\vert >r_{n}v_{n}^{-1/2}$.

Therefore,
\begin{eqnarray*}
\vert{-}z^{-2}+v_{n}\vert
&=&v_{n}\vert z\vert ^{-2}\vert z-v_{n}^{-1/2}\vert\vert z+v_{n}^{-1/2}\vert
\\
&>& v_{n}v_{n}r_{n}v_{n}^{-1/2}=r_{n}v_{n}^{3/2}.
\end{eqnarray*}

Since $r_{n}+v_{n}^{-1/2}<2v_{n}^{-1/2}$, assumption (i) implies
that $r_{n}+v_{n}^{-1/2}<m_{n}/2$. Hence, the circle $\vert
z\vert=v_{n}^{-1/2}+r_{n}$ lies entirely in the area where
formula (%
\ref{K_prime_estimate}) applies to $K_{n}^{\prime}( z) $.
Consequently, using (\ref{K_prime_estimate}), we can estimate
\[
\vert K_{n}^{\prime}( z) -( -z^{-2}+v_{n})
\vert\leq4D_{n}v_{n}^{-1/2},
\]
where we used the fact that
$\vert z\vert\leq2v_{n}^{-1/2}$. By assumption (ii),
$r_{n}\geq4D_{n}v_{n}^{-2}$, therefore $v_{n}^{3/2}r_{n}\geq
4D_{n}v_{n}^{-1/2}$ and Rouch\'e's theorem is applicable
to the pair of $K_n^{\prime}(z)$ and $-z^{-2}+v_n$. Both
$K_{n}^{\prime}( z) $ and $-z^{-2}+v_{n}$ have only one pole,
of order two, in $\vert z\vert\leq v^{-1/2}-r_{n}$
and the function $-z^{-2}+v_{n}$ has no zeros inside
$\vert z\vert\leq v^{-1/2}-r_{n}$. Therefore, Rouch\'e's
theorem implies that there are no zeros
of $K_{n}^{\prime}( z) $ inside $\vert z\vert \leq v^{-1/2}-r_{n}$.
(Rouch\'e's theorem is often formulated only for holomorphic
functions, but as a consequence of the argument principle (see, e.g., Theorems
II.2.3 and II.2.4 in \cite{markushevich77}), it can be easily reformulated
for meromorphic functions. In this form, it claims that a meromorphic
function, $f(z)$, has the same difference between the number of zeros and
number of poles inside a curve $\gamma$ as another meromorphic
function, $g(z)$, provided that $\vert f( z) \vert>\vert g( z) -f( z) \vert$.
For this formulation see, e.g., \cite{hille62}, Theorem 9.2.3.)
\end{pf*}
\begin{condition}
Assume, in the following, that $r_{n}=4D_{n}/v_{n}^{2}$.
\end{condition}

We now use our knowledge about the location of critical points of
$K_{n}( z) $ to investigate how it behaves on the real interval
around zero.
\begin{lemma}\label{real_map}
Suppose that $m_{n}>4/\sqrt{v_{n}}$ and
$D_{n}/v_{n}^{3/2}\leq
1/8 $. Then $K_{n}( z) $ maps the set $[ -1/\sqrt{v_{n}}
+r_{n},0) \cup( 0,1/\sqrt{v_{n}}-r_{n}] $ in a one-to-one
fashion onto a set that contains the union of two intervals
$(-\infty,-2\sqrt{v_{n}}-cD_{n}/v_{n}) \cup( 2\sqrt{v_{n}}+cD_{n}/v_{n},\infty) $,
where $c$ is a constant that does not depend on $n$.
\end{lemma}
\begin{rem*}
For example, $c=5$ will work.
\end{rem*}
\begin{pf*}{Proof of Lemma \ref{real_map}}
The assumption that $m_{n}>4/\sqrt{v_{n}}$ ensures that the
power series for $K_{n}( z) $ converges in $\vert
z\vert\leq4/\sqrt{v_{n}}$, $z\neq0$. Note that
$K_{n}( z) $ is real-valued on the set $I=[ -1/\sqrt{v_{n}}+r_{n},0)
\cup( 0,1/\sqrt{v_{n}}-r_{n}] $ because this
set belongs to the area where the series for $K_{n}( z) $
converges and the coefficients of this series are real. Moreover, by
Lemma \ref{K_zeros}, there are no critical points of $K_{n}( z) $
on~$I$ [i.e., for every $z\in I$, $K_{n}^{\prime}( z) \neq0$],
therefore $K_{n}( z) $ must be strictly monotonic on subintervals
$[ -1/\sqrt{v_{n}}+r_{n},0) $ and $( 0,1/\sqrt{v_{n}}-r_{n}] $.
Consequently, $K_{n}( I) =( -\infty,K_{n}(-1/\sqrt{v_{n}}+r_{n})]
\cup[ K_{n}( 1/\sqrt{v_{n}}-r_{n}) ,\infty) $. We claim that
$K_{n}( 1/\sqrt{v_{n}}-r_{n}) \leq2\sqrt{v_{n}}+5D_{n}/v_{n}$ and
$K_{n}( -1/\sqrt{v_{n}}+r_{n}) \geq-2\sqrt{v_{n}}-5D_{n}/v_{n}$.

Indeed, if we write
\[
K_{n}( z) =\frac{1}{z}+v_{n}z+h( z) ,
\]
then
\[
K_{n}\biggl( \frac{1}{\sqrt{v_{n}}}-r_{n}\biggr)
= \sqrt{v_{n}}\frac{1}{1-r_{n}\sqrt{v_{n}}}
+ \sqrt{v_{n}}\bigl( 1-r_{n}\sqrt{v_{n}}\bigr)
+ h\biggl(\frac{1}{\sqrt{v_{n}}}-r_{n}\biggr) .
\]
According to our assumption, $r_{n}\sqrt{v_{n}}=4D_{n}/v_{n}^{3/2}<1/2$.
Therefore, we can estimate
\[
\frac{1}{1-r_{n}\sqrt{v_{n}}}\leq1+2r_{n}\sqrt{v_{n}}
\]
and
\[
K_{n}\biggl( \frac{1}{\sqrt{v_{n}}}-r_{n}\biggr)
\leq2\sqrt{v_{n}} +r_{n}v_{n}
+ \biggl\vert h\biggl( \frac{1}{\sqrt{v_{n}}}-r_{n}\biggr)\biggr\vert.
\]
We can estimate the last term using Lemma \ref{Lemma_K_estimate} as
\[
h\biggl( \frac{1}{\sqrt{v_{n}}}-r_{n}\biggr)
\leq D_{n} \biggl\vert\frac{1}{\sqrt{v_{n}}}\biggr\vert^{2}
= D_{n}/v_{n}.
\]
Combining all of this and substituting $r_{n}=4D_{n}/v_{n}^{2}$, we get
\[
K_{n}\biggl( \frac{1}{\sqrt{v_{n}}}-r_{n}\biggr)
\leq2\sqrt{v_{n}} + 5D_{n}/v_{n}.
\]
Similarly, we can derive that
\[
K_{n}\biggl( -\frac{1}{\sqrt{v_{n}}}+r_{n}\biggr)
\geq-2\sqrt{v_{n}}-5D_{n}/v_{n}.
\]\upqed
\end{pf*}

From the previous lemma, we can conclude that $K_{n}( z) $
has a differentiable inverse defined on $( -\infty,-2\sqrt{v_{n}}
-cD_{n}/v_{n}) \cup( 2\sqrt{v_{n}}+cD_{n}/v_{n},\infty) $.
We can extend this conclusion to an open complex neighborhood of this
interval. This is achieved in the next two lemmas.
\begin{lemma}
As in the previous lemma, suppose that $m_{n}>4/\sqrt{v_{n}}$ and
$D_{n}/v_{n}^{3/2}\leq1/8$. Let $z$ be an arbitrary point of the
interval $%
[ -1/\sqrt{v_{n}}+r_{n},1/\sqrt{v_{n}}-r_{n}] $. Then we can find
a neighborhood $U_{z}$ of $z$ and a neighborhood $W_{w}$ of
$w=K_{n}(z) $ such that $K_{n}$ is a one-to-one map of $U_{z}$ onto
$W_{w}$ and the inverse map~$K_{n}^{-1}$ is holomorphic everywhere in $W_{w}$.
\end{lemma}
\begin{pf}
Since the power series for $K_{n}( z) -z^{-1}$ converges
in $\vert z\vert\leq4/\sqrt{v_{n}}$, the function
$K_{n}(z) $ is holomorphic in $\vert z\vert\leq4/\sqrt{v_{n}}$,
$z\neq0$. In addition, by Lemma~\ref{K_zeros}, $z\in[ -1/\sqrt
{v_{n}}+r_{n},1/\sqrt{v_{n}}-r_{n}] $ is not a critical point of
$K_{n}(z) $. Therefore, for $z\neq0$, the conclusion of the lemma follows
from Theorems II.3.1 and II.3.2 in~\cite{markushevich77}. For $z=0$, the
argument is parallel to the argument in Markushevich, except for a different
choice of local coordinates. Indeed, $f( z) =1/K_{n}(z) $ is
holomorphic at $z=0$, it maps $z=0$ to $w=0$ and
$f^{\prime}( z) =1\neq0$ at $z=0$. Therefore, Theorems II.3.1 and II.3.2
in \cite{markushevich77} are applicable to $f( z) $ and it
has a well-defined holomorphic inverse in a neighborhood of $w=0$. This implies
that $K_{n}(z)$ has a well-defined holomorphic inverse in a
neighborhood of $\infty$, given by the formula $K_{n}^{-1}( z) =f^{-1}(1/z) $.
\end{pf}
\begin{lemma}\label{Global_inverse}
Local inverse $K_{n}^{-1}( z) $ defined in
the previous lemma is a restriction of a function $G_{n}( z) $
which is defined and holomorphic everywhere in a neighborhood of
$I=\{\infty\} \cup( -\infty,-2v_{n}^{1/2}-cD_{n}/v_{n}]
\cup[ 2v_{n}^{1/2}+cD_{n}/v_{n},\infty) $. The function
$G_{n}(z) $ is the inverse of $K_{n}( z) $ in this neighborhood.
\end{lemma}
\begin{pf}
By Lemma \ref{real_map}, for every point $w\in I$, we can
find a unique $z\in[ -1/\sqrt{v_{n}}+r_{n},1/\sqrt{v_{n}}-r_{n}]$
such that $K_{n}( z) =w$. Let $U_{z}$ and $W_{w}$ be the neighborhoods
defined in the previous lemma. Also, let us write
$(K_{n}^{-1},W_{w}) $ to denote the local inverses defined in the
previous lemma together with their areas of definition. Our task is to prove
that these local inverses can be joined to form an analytic function
$K_{n}^{-1}$, well defined everywhere in a neighborhood of $I$. We will
do this in several steps.

First, an examination of the proof of the previous lemma and Theorem II.3.1
in~\cite{markushevich77} shows that we can take each $U_{z}$ in the
form of a disc. Let $\widetilde{U}_{z}=U_{z}/3$, that is, define
$\widetilde{U}_{z}$ as a disc that has the same center, but
radius one third that of~$U_{z}$. Define
$\widetilde{W}_{w}$ as $K_{n}( \widetilde{U}_{z})$.
These new sets are more convenient because of the
property that if
$\widetilde{U}_{z_{1}}\cap\widetilde{U}_{z_{2}}\neq\varnothing$,
then either $\widetilde{U}_{z_{1}}\cup\widetilde{U}_{z_{2}}\subset U_{z_{1}}$
or $\widetilde{U}_{z_{1}}\cup\widetilde{U}_{z_{2}}\subset U_{z_{2}}$. In
particular, this means that if $\widetilde{U}_{z_{1}}\cap\widetilde{U}
_{z_{2}}\neq\varnothing$, then $K_{n}( z) $ is a one-to-one
map of $\widetilde{U}_{z_{1}}\cup\widetilde{U}_{z_{2}}$ onto
$\widetilde{W}_{w_{1}}\cup\widetilde{W}_{w_{2}}$. This is convenient because $K_{n}$ is
one-to-one not only on a particular neighborhood
$\widetilde{U}_{z_{1}}$, but
also on the union of every two intersecting neighborhoods
$\widetilde{U}_{z_{1}}$ and $\widetilde{U}_{z_{2}}$.
Let us call this the \emph{extended invertibility} property.

Next, define an even smaller $\widetilde{\widetilde{U}}_{z}$ with the following
properties: (1) $\widetilde{\widetilde{U}}_{z}\subset\widetilde{U}_{z}$;
(2)~$\widetilde{\widetilde{W}}_{w}=:K_{n}( \widetilde{\widetilde{U}}_{z}) $
is either an open disc for $z\neq0$ or the set
$\vert w\vert>R$ for $z=0$; and (3) $0\notin\widetilde{\widetilde{W}}_{w}$.
This is easy to achieve by taking an appropriate open subset of
$\widetilde{W}_{w}$ as $\widetilde{\widetilde{W}}_{w}$ and applying $K_{n}^{-1}$. Note
that the property of the previous paragraph
(i.e., extended invertibility) remains valid for the new
sets $\widetilde{\widetilde{U}}_{z}$.

Discs $\widetilde{\widetilde{W}}_{w}$ form an open cover of $I$ and the
corresponding sets $\widetilde{\widetilde{U}}_{z}$ form an open cover
for $K_{n}^{-1}( I) $, which is a closed interval contained in
$[-1/\sqrt{v_{n}}+r_{n},1/\sqrt{v_{n}}-r_{n}] $. Let $U_{_{i}}$,
$i=0,\ldots,N$, be a finite cover of $K_{n}^{-1}( I) $ selected
from $\{ \widetilde{\widetilde{U}}_{z}\} $. [Recall that $K_{n}^{-1}$
is well defined on the interval $I$ by Lemma \ref{real_map}, hence
$K_{n}^{-1}( I) $ is well defined.] We can find a finite cover
due to the compactness of $K_{n}^{-1}( I) $. Further, let
$W_{_{i}}=:K_{n}( U_{i}) $ be the corresponding cover of $I$,
selected from $\{ \widetilde{\widetilde{W}}_{z}\} $. For
convenience, let $W_{0}$ denote the set
$\widetilde{\widetilde{W}}_{w}$ for $w=\infty$. Finally,
let $R=\bigcup_{i=0}^{N}U_{i}$ and $S=\bigcup_{i=0}^{N}W_{i}$.
Sets $R$ and $S$ are illustrated in Figure \ref{f1}.

\begin{figure}

\includegraphics{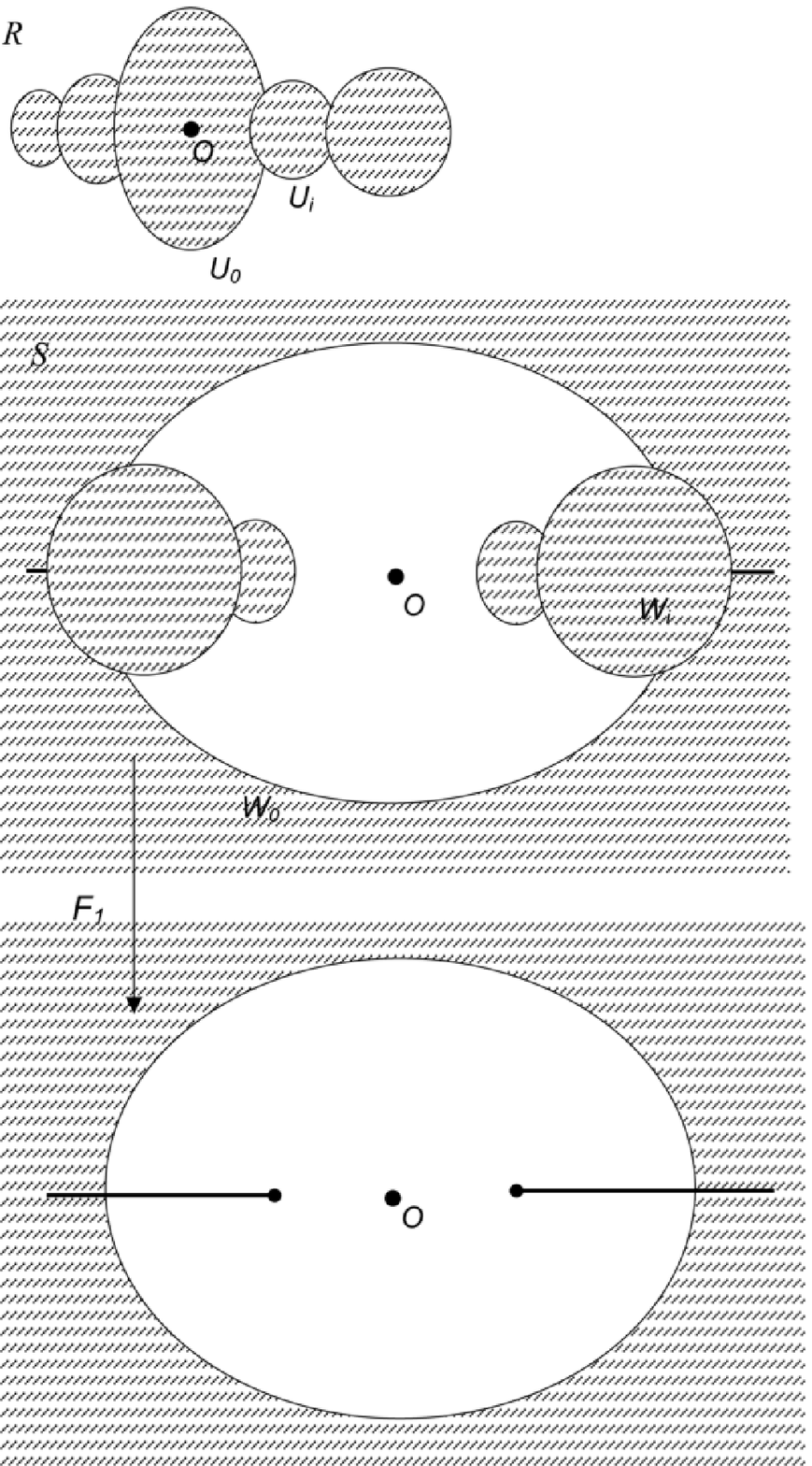}

  \caption{\label{f1}}
\end{figure}

Clearly, $S$ is open. We aim to prove that $S$ is simply connected in the
extended complex plane $\mathbb{C}\cup\{ \infty\} $. For this
purpose, let us define the deformation retraction~$F_{1}$ of the set $S$
as  follows: (1) if $z\in W_{0}$, then $z\rightarrow z$; (2)
if $z\notin W_{0}$, then $z\rightarrow\operatorname{Re}z+( 1-t) \operatorname{Im}z$.
Here, parameter $t$ changes from 0 to 1.
(For the definition and properties of deformation
retractions, see, e.g., \cite{hatcher02}; the definition is on page 2 and
the main property is in Proposition 1.17.) This retraction reduces $S$
to a
homotopically equivalent set $S^{\prime}$ that consists of $W_{0}$ and two
intervals of the real axis that do not include $0$. We can then use
another deformation retraction $F_{2}$ that sends $z$ to $(1-t)^{-1}z$. This
retraction reduces $S^{\prime}$ to $S^{\prime\prime}=\{ \infty\} $,
which is evidently simply connected.

We know that there is a holomorphic inverse $K_{n}^{-1}(z)$ defined on each
of~$W_{i}$. Starting from one of these domains, say $W_{0}$, we can
analytically continue $K_{n}^{-1}( z) $ to every other $W_{i}$.
Indeed, take points $z_{0}\in U_{0}$ and $z_{i}\in U_{i}$ and connect them
by a path that lies entirely in $R=\bigcup_{i=0}^{N}U_{i}$. This path
corresponds to a chain~$\{ U_{k_{s}}\} $, $s=1,\ldots,n$, that
connects $U_{0}$ and $U_{i}$. That is, $U_{k_{1}}=U_{0}$, $U_{k_{n}}=U_{i}$
and \mbox{$U_{k_{j}}\cap U_{k_{j+1}}\neq\varnothing$}. The corresponding
$W_{k_{s}}=K_{n}( U_{k_{s}}) $ form a chain that connects $W_{0}$
and $W_{j}$, that is, $W_{k_{1}}=W_{0}$, $W_{k_{n}}=W_{i}$
and $W_{k_{j}}\cap W_{k_{j+1}}\neq\varnothing$. By our construction, this chain
of sets $W_{k_{s}}$ has the property that $K_{n}^{-1}(W_{k_{j}})
\cap K_{n}^{-1}( W_{k_{j+1}}) =U_{k_{j}}\cap U_{k_{j+1}}\neq\varnothing$.

Consider two adjacent sets, $W_{k_{j}}$ and $W_{k_{j+1}}$, in this chain.
Then the corresponding local inverses $( K_{n}^{-1},W_{k_{j}}) $
and $( K_{n}^{-1},W_{k_{j+1}}) $, which were defined in the
previous lemma, coincide on an open nonempty set. Indeed,
$K_{n}(U_{k_{j}}\cap U_{k_{j+1}}) \subset K_{n}( U_{k_{j}})
\cap K_{n}( U_{k_{j+1}}) =W_{k_{j}}\cap W_{k_{j+1}}$, therefore the
functions $( K_{n}^{-1},W_{k_{j}}) $ and
$(K_{n}^{-1},W_{k_{j+1}}) $ are both well defined on
$K_{n}(U_{k_{j}}\cap U_{k_{j+1}}) $. Moreover, they must coincide on
$K_{n}( U_{k_{j}}\cap U_{k_{j+1}}) $. Indeed, by construction,
$U_{k_{j}}\cap U_{k_{j+1}}\neq\varnothing$ and therefore, by the extended
invertibility property, $K_{n}$ is one-to-one on
$U_{k_{j}}\cup U_{k_{j+1}}$.
Hence, there cannot exist two different $z$ and
$z^{\prime}$ in $U_{k_{j}}\cup U_{k_{j+1}}$
that would map to the same point in
$K_{n}(U_{k_{j}}\cap U_{k_{j+1}}) $. Hence,
$( K_{n}^{-1},W_{k_{j}}) $ and
$( K_{n}^{-1},W_{k_{j+1}}) $ must coincide on
$K_{n}(U_{k_{j}}\cap U_{k_{j+1}}) $, which is open and nonempty.

Using the property that if two analytical functions coincide on an open
set,
then each of them is an analytic continuation of the other, we conclude that
the local inverse $( K_{n}^{-1},W_{k_{j}}) $ can be analytically
continued to $W_{k_{j+1}}$, where it coincides with the local inverse
$( K_{n}^{-1},W_{k_{j+1}}) $. Therefore, at least one analytic
continuation of $( W_{0},K_{n}^{-1}) $ is well defined everywhere
on $S$ and has the property that when restricted to each of $W_{j}$, it
coincides with a local inverse of $K_{n}( z) $ defined in the
previous lemma. Since $S$ is simply connected, the analytic
continuation is
unique, that is, it does not depend on the choice of the chain of the
neighborhoods that connect $W_{0}$ and $W_{j}$.

Let us denote the function resulting from this analytic continuation as
$G_{n}( z) $. By construction, it is unambiguously defined for
every $W_{j}$ and the restrictions of $G_{n}( z) $ to $W_{j}$
coincide with $( K_{n}^{-1},W_{j}) $. Therefore,
$G_{n}(z) $ satisfies the relations $K_{n}( G_{n}( z))=z$ and
$G_{n}( K_{n}( z) ) =z$ everywhere on
$R=\bigcup_{i=0}^{N}U_{i}$ and on
$S=\bigcup_{i=0}^{N}W_{i}$. Since $S$ is an open
neighborhood of $I$, every claim of the lemma is proved.
\end{pf}
\begin{lemma}\label{Global_is_Cauchy}
The function $G_{n}( z) $ constructed in
the previous lemma is the Cauchy transform of $S_{n}$.
\end{lemma}

By construction, $G_{n}^{-1}( z) $ is the inverse of
$K_{n}(z) $ in a neighborhood of $\{ \infty\} \cup(-\infty,-2v_{n}^{1/2}-cD_{n}/v_{n})
\cup(2v_{n}^{1/2}+cD_{n}/v_{n},\infty) $. In particular, it is the
inverse of $K_{n}( z) $ in a neighborhood of infinity. Therefore, in this
neighborhood, it has the same power expansion as the Cauchy transform of
$S_{n}$. Therefore, it coincides with the Cauchy transform of $S_{n}$ in this
neighborhood. Next, we apply the principle that if two analytical functions
coincide in an open domain, then they coincide at every point where they can
be continued analytically.

It remains to apply Lemma \ref{G_boundedness_lemma} in order to
obtain the following theorem.
\begin{theorem}\label{main_theorem copy(1)}
Suppose that \textup{(i)} $\lim\inf m_{n}\sqrt{v_{n}}>4$ and
\textup{(ii)} $\lim\sup_{n\rightarrow\infty}D_{n}/\break v_{n}^{3/2}\leq1/8$.
Then for all sufficiently large $n$, the support of $\mu_{n}$ belongs to
\[
I=\bigl( -2\sqrt{v_{n}}-cD_{n}/v_{n},2\sqrt{v_{n}}+cD_{n}/v_{n}\bigr) ,
\]
where $c>0$ is an absolute constant (e.g., $c=5$).
\end{theorem}
\begin{pf*}{Proof of Theorem \ref{main_theorem copy(1)}}
Let us collect the facts that we know about
the function $G_{n}( z) $ defined in Lemma
\ref{Global_inverse}. First, by Lemma \ref{Global_is_Cauchy}, it is the
Cauchy transform of a bounded random variable, $S_{n}$. Second, by Lemma
\ref{Global_inverse}, it is holomorphic at all $z\in\mathbb{R}$ such that
$\vert z\vert>2v_{n}^{1/2}+cD_{n}/v_{n}$. Using Lemma \ref{G_boundedness_lemma},
we conclude that the distribution of $S_{n}$ is
supported on the interval $[-2v_{n}^{1/2}-cD_{n}/v_{n},2v_{n}^{1/2}+cD_{n}] $.
\end{pf*}

If we take $R_{n,i}=2L_{n,i}$ and $m_{n,i}=( 4L_{n,i}) ^{-1}$,
then assumption (i) is equivalent to
\[
\mathop{\lim\inf}_{n\rightarrow\infty}\min_{i}\frac{\sqrt{v_{n}}}{4L_{n,i}}>4,
\]
which is equivalent to
\[
\mathop{\lim\sup}_{n\rightarrow\infty}\frac{L_{n}}{\sqrt{v_{n}}}<16.
\]
From (ii), we obtain
\[
1/8\geq\mathop{\lim\sup}_{n\rightarrow\infty}\frac{\sum_{i=1}^{k_{n}}R_{n,i}%
( m_{n,i}) ^{-2}}{v_{n}^{3/2}}=\mathop{\lim\sup}_{n\rightarrow
\infty}%
\frac{32\sum_{i=1}^{k_{n}}L_{n,i}^{3}}{v_{n}^{3/2}},
\]
which is equivalent to
\[
\mathop{\lim\sup}_{n\rightarrow\infty}\frac{T_{n}}{v_{n}^{3/2}}\leq1/256.
\]

Finally, note that the condition $\lim\sup_{n\rightarrow\infty}
T_{n}/v_{n}^{3/2}\leq2^{-12}$ implies that
$\lim\sup_{n\rightarrow\infty}L_{n}/\sqrt{v_{n}}<16$. Therefore, Theorem
\ref{main_theorem} is a consequence of Theorem~\ref{main_theorem copy(1)}.

\section*{Acknowledgment}

I would like to express my gratitude to Diana Bloom for her editorial help.\vadjust{\goodbreak}

\printaddresses

\end{document}